\magnification 1200

%
%
\newdimen\FigSize       \FigSize=.9\hsize 
%
\newskip\abovefigskip   \newskip\belowfigskip
\gdef\epsfig#1;#2;{\par\vskip\abovefigskip\penalty -500
   {\everypar={}\epsfxsize=#1\nd \centerline{\epsfbox{#2}}}%
    \vskip\belowfigskip}%
%
\newskip\figtitleskip
\gdef\tepsfig#1;#2;#3{\par\vskip\abovefigskip\penalty -500
   {\everypar={}\epsfxsize=#1\nd
    \vbox
      {\centerline{\epsfbox{#2}}\vskip\figtitleskip
       \centerline{\figtitlefont#3}}}%
    \vskip\belowfigskip}%
%
\newcount\FigNr \global\FigNr=0
\gdef\nepsfig#1;#2;#3{\global\advance\FigNr by 1
   \tepsfig#1;#2;{Figure\space\the\FigNr.\space#3}}%
%
%
%
\gdef\ipsfig#1;#2;{
   \midinsert{\everypar={}\epsfxsize=#1\nd
              \centerline{\epsfbox{#2}}}
   \endinsert}%
%
\gdef\tipsfig#1;#2;#3{\midinsert
   {\everypar={}\epsfxsize=#1\nd
    \vbox{\centerline{\epsfbox{#2}}%
          \vskip\figtitleskip
          \centerline{\figtitlefont#3}}}\endinsert}%
%
\gdef\nipsfig#1;#2;#3{\global\advance\FigNr by1%
  \tipsfig#1;#2;{Figure\space\the\FigNr.\space#3}}%
\newread\epsffilein    
\newif\ifepsffileok    
\newif\ifepsfbbfound   
\newif\ifepsfverbose   
\newdimen\epsfxsize    
\newdimen\epsfysize    
\newdimen\epsftsize    
\newdimen\epsfrsize    
\newdimen\epsftmp      
\newdimen\pspoints     
\pspoints=1bp          
\epsfxsize=0pt         
\epsfysize=0pt         
\def\epsfbox#1{\global\def\epsfllx{72}\global\def\epsflly{72}%
   \global\def\epsfurx{540}\global\def\epsfury{720}%
   \def\lbracket{[}\def\testit{#1}\ifx\testit\lbracket
   \let\next=\epsfgetlitbb\else\let\next=\epsfnormal\fi\next{#1}}%
\def\epsfgetlitbb#1#2 #3 #4 #5]#6{\epsfgrab #2 #3 #4 #5 .\\%
   \epsfsetgraph{#6}}%
\def\epsfnormal#1{\epsfgetbb{#1}\epsfsetgraph{#1}}%
\def\epsfgetbb#1{%
%
%
\openin\epsffilein=#1
\ifeof\epsffilein\errmessage{I couldn't open #1, will ignore it}\else
%
%
   {\epsffileoktrue \chardef\other=12
    \def\do##1{\catcode`##1=\other}\dospecials \catcode`\ =10
    \loop
       \read\epsffilein to \epsffileline
       \ifeof\epsffilein\epsffileokfalse\else
%
%
          \expandafter\epsfaux\epsffileline:. \\%
       \fi
   \ifepsffileok\repeat
   \ifepsfbbfound\else
    \ifepsfverbose\message{No bounding box comment in #1; using
defaults}\fi\fi
   }\closein\epsffilein\fi}%
%
%
\def\epsfsetgraph#1{%
   \epsfrsize=\epsfury\pspoints
   \advance\epsfrsize by-\epsflly\pspoints
   \epsftsize=\epsfurx\pspoints
   \advance\epsftsize by-\epsfllx\pspoints
%
%
   \epsfxsize\epsfsize\epsftsize\epsfrsize
   \ifnum\epsfxsize=0 \ifnum\epsfysize=0
      \epsfxsize=\epsftsize \epsfysize=\epsfrsize
%
arithmetic!
%
     \else\epsftmp=\epsftsize \divide\epsftmp\epsfrsize
       \epsfxsize=\epsfysize \multiply\epsfxsize\epsftmp
       \multiply\epsftmp\epsfrsize \advance\epsftsize-\epsftmp
       \epsftmp=\epsfysize
       \loop \advance\epsftsize\epsftsize \divide\epsftmp 2
       \ifnum\epsftmp>0
          \ifnum\epsftsize<\epsfrsize\else
             \advance\epsftsize-\epsfrsize \advance\epsfxsize\epsftmp
\fi
       \repeat
     \fi
   \else\epsftmp=\epsfrsize \divide\epsftmp\epsftsize
     \epsfysize=\epsfxsize \multiply\epsfysize\epsftmp
     \multiply\epsftmp\epsftsize \advance\epsfrsize-\epsftmp
     \epsftmp=\epsfxsize
     \loop \advance\epsfrsize\epsfrsize \divide\epsftmp 2
     \ifnum\epsftmp>0
        \ifnum\epsfrsize<\epsftsize\else
           \advance\epsfrsize-\epsftsize \advance\epsfysize\epsftmp \fi
     \repeat
   \fi
%
%
   \ifepsfverbose\message{#1: width=\the\epsfxsize,
height=\the\epsfysize}\fi
   \epsftmp=10\epsfxsize \divide\epsftmp\pspoints
   \vbox to\epsfysize{\vfil\hbox to\epsfxsize{%
      \includegraphics{#1}%
      \hfil}}%
\epsfxsize=0pt\epsfysize=0pt}%
%
%
{\catcode`\%=12
\global\let\epsfpercent=
%
%
\long\def\epsfaux#1#2:#3\\{\ifx#1\epsfpercent
   \def\testit{#2}\ifx\testit\epsfbblit
      \epsfgrab #3 . . . \\%
      \epsffileokfalse
      \global\epsfbbfoundtrue
   \fi\else\ifx#1\par\else\epsffileokfalse\fi\fi}%
%
%
\def\epsfgrab #1 #2 #3 #4 #5\\{%
   \global\def\epsfllx{#1}\ifx\epsfllx\empty
      \epsfgrab #2 #3 #4 #5 .\\\else
   \global\def\epsflly{#2}%
   \global\def\epsfurx{#3}\global\def\epsfury{#4}\fi}%
%
%
\def\epsfsize#1#2{\epsfxsize}%
%
%

\epsfverbosetrue                        
\abovefigskip=\baselineskip             
\belowfigskip=\baselineskip             
\global\let\figtitlefont\bf             
\global\figtitleskip=.5\baselineskip    

\newread\rfaux
\openin\rfaux=\jobname.AUX
\ifeof\rfaux \message{Can't find \jobname.AUX!!!}
   \else \closein\rfaux
     \input \jobname.AUX \fi

\newwrite\wfaux
\immediate\openout\wfaux=\jobname.AUX

\newcount\eqnum  \newif\ifeqnerr

\def\neweq#1{{\global\advance\eqnum by 1
    \edef\chqtmp{(\the\eqnum)}
    \immediate\write\wfaux{\string\def \string#1{\chqtmp}}
    \checkeq#1}}

\def\checkeq#1{\ifx#1\chqtmp\global\eqnerrfalse
    \else \global\eqnerrtrue
       \ifx#1\zundefined
           \def\zz{NOT DEFINED}\global\let#1\chqtmp
       \else \def\zz{preassigned #1}\fi
        \message{ALLOCATION ERROR for \string#1:
             \zz, should be \chqtmp.}\fi}

\def\chq#1{\neweq#1\ifeqnerr***\fi#1}
\def\cheqno{\eqno\chq}

\def\pagemac#1{
    \write\wfaux{\string\def \string#1{\the\pageno}}
    \ifx#1\zundefined
      \message{Page macro \string#1 NOT DEFINED.}\gdef#1{***}\fi}


\def\neweqx#1#2{{\edef\chqtmp{(\the\eqnum#1)}
    \immediate\write\wfaux{\string\def \string#2{\chqtmp}}
    \checkeq#2}}

\font\tenmsb=msbm10   
\font\sevenmsb=msbm7
\font\fivemsb=msbm5
\newfam\msbfam
\textfont\msbfam=\tenmsb
\scriptfont\msbfam=\sevenmsb
\scriptscriptfont\msbfam=\fivemsb
\def\Bbb#1{\fam\msbfam\relax#1}
\let\nd\noindent 
\def\qed{\hbox{\hskip 6pt\vrule width6pt height7pt depth1pt \hskip1pt}}
\def\natural{{\rm I\kern-.18em N}}
\def\a{{\bf a}}
\def\b{{\bf b}}
\def\s{{\bf s}}
\def\t{{\bf t}}
\def\A{{\cal A}}
\def\B{{\cal B}}

\def\E{{\Bbb E}}
\def\G{{\cal G}}

\def\L{{\cal L}}
\def\N{{\Bbb N}}
\def\P{{\cal P}}

\def\R{{\Bbb R}}

\def\Z{{\Bbb Z}}
\def\RO{{\rm RO}}
\def\chix{{\raise.5ex\hbox{$\chi$}}}
\def\chixa{{\chix\lower.2em\hbox{$_A$}}}

\def\real{{\rm I\kern-.2em R}}
\def\integer{{\rm Z\kern-.32em Z}}
\def\complex{\kern.1em{\raise.47ex\hbox{
            $\scriptscriptstyle |$}}\kern-.40em{\rm C}}
\def\vs#1 {\vskip#1truein}
\def\hs#1 {\hskip#1truein}

\def\Month{\ifcase\number\month \relax\or January \or February \or
  March \or April \or May \or June \or July \or August \or September
  \or October \or November \or December \else \relax\fi }
\def\date{\Month \the\day, \the\year}

  \hsize=6truein        \hoffset=.25truein 
  \vsize=8.8truein      
  \pageno=1     \baselineskip=12pt
  \parskip=0 pt         \parindent=20pt
  \overfullrule=0pt     \lineskip=0pt   \lineskiplimit=0pt
  \hbadness=10000 \vbadness=10000 
\nd
\pageno=0

\footline{\ifnum\pageno=0\hss\else\hss\tenrm\folio\hss\fi}
\hbox{}
\vskip 1truein\centerline{{\bf CONJUGACIES FOR TILING DYNAMICAL SYSTEMS}}
\vskip .2truein\centerline{by}
\vskip .2truein
\centerline{
Charles Holton
\footnote{*}{Research supported in part by NSF Vigre Grant DMS-0091946},
Charles Radin
\footnote{**}{Research supported in part by NSF Grant DMS-0071643 and
\hfill\break \indent Texas ARP Grant 003658-158\hfil},
and Lorenzo Sadun}
\vskip .2truein
\centerline{Department of Mathematics, University of Texas, Austin, TX}

\vs1.5
\centerline{{\bf Abstract}}
\vs.2 \nd
We consider tiling dynamical systems and topological conjugacies
between them. We prove that
the criterion of being finite type is invariant under topological
conjugacy. For substitution tiling systems under rather general
conditions, including the Penrose and pinwheel systems, we show that
substitutions are invertible and that
conjugacies are generalized sliding block codes. 
\vs1 \nd
\vs.2
\vfill\eject

\nd
{\bf I.\ Notation and main results}
\vs.1

We begin with a definition of tiling dynamical systems, in sufficient
generality for this work. Let $\A$ be a nonempty finite collection of
compact connected sets in the Euclidean space $\E^d$, sets with dense
interior and boundary of zero volume. Let $X(\A)$ be the set of all
tilings of $\E^d$ by congruent copies, which we call tiles, of the
elements of the ``alphabet'' $\A$.  We assume $X(\A)$ is nonempty,
which is automatic for the special class of substitution tiling
systems on which we will concentrate below. We label the ``types'' of
tiles by the elements of $\A$.  We endow $ X(\A)$ with the metric
$$
m[x, y]\equiv \sup_{n\ge 1}{1\over n}m_{H}[B_{n}\cap \partial
x,B_{n}\cap \partial y], \cheqno\ea
$$
where ${B}_{n}$ denotes the open ball of radius $n$ centered at the
origin {\bf O} of $\E^d$, and $\partial x$ the union of the boundaries
of all tiles in $x$. (A ball centered at $\a$ is denoted $B_n(\a)$.)
The Hausdorff metric $m_{H}$ is defined as follows. Given two compact
subsets $P$ and $Q$ of $\E^d$, $m_{H}[P,Q] = \max \{ {\tilde m} (P,Q),
{\tilde m} (Q,P)\}$, where
$${\tilde m} (P,Q) =  \sup_{p \in P} \inf_{q \in Q} \|p - q\|, \cheqno\eb
$$
with $\|w\|$ denoting the usual Euclidean norm of $w$.

Under the metric $m$ two tilings are close if they agree, up to a
small Euclidean motion, on a large ball centered at the origin.  The
converse is also true for tiling systems with finite local complexity
(as defined below): closeness implies agreement, up to small Euclidean
motion, on a large ball centered at the origin [see RaS1].  Although
the metric $m$ depends on the location of the origin, the topology
induced by $m$ is Euclidean invariant.  A sequence of tilings
converges in the metric $m$ if and only if its restriction to every
compact subset of $\E^d$ converges in $m_H$.  It is not hard to show
[RW] that $X(\A)$ is compact and that the natural action of the
connected Euclidean group $\G_E$ on $X(\A)$, $(g,x)\in \G_E\times
X(\A)\mapsto g[x]\in X(\A)$, is continuous.

To include certain examples it is useful to generalize the above
setup, to use what is sometimes called ``colored tiles''. To make the
generalization we assign a ``color'' from some finite set to each
element of $\A$, represented on each tile by a ``color marking'', a
line segment in the interior of the tile, of different length for
different colors. We then redefine $\partial x$ as the union of the
tile boundaries and color markings in the tiling $x$.
\vs.1 \nd
{\bf Definition 1.} A tiling dynamical system is the action of $\G_E$ on
a closed,
\break \nd $\G_E$-invariant subset of $X(\A)$.
\vs.1
We emphasize the close connection between such dynamical systems and
subshifts. A subshift with $\Z^d$-action is the natural action of
$\Z^d$ on a compact, $\Z^d$-invariant subset $X$ of $\B^{\Z^d}$, for some
nonempty finite set $\B$. If we
associate with each element of $\B$ a ``colored'' unit cube in $\E^d$,
the face-to-face tilings of $\E^d$ by those arrays of such cubes
corresponding to the subshift $X$ gives a tiling dynamical system
which is basically the suspension of the subshift $X$ (but
with rotations of the entire tiling also permitted).

A significant difference between subshifts and tiling dynamical
systems is that for (nontrivial) subshifts the group acts on a Cantor
set, while the space is typically connected for interesting tiling
systems. In fact, the spaces for different tiling systems need not be
homeomorphic.

A major objective in dynamics is the classification of interesting
subclasses up to topological conjugacy. For the class of subshifts a
central theorem, due to Curtis, Lyndon and Hedlund, shows that a
topological conjugacy can be represented by a sliding block code (see
[LM]). For tiling dynamical systems there is a natural analogue of
such a representation for which we use the same term.
(Such maps are called ``local'' in [P] and are closely related to
mutual local derivability [BSJ].)
\vs.1 \nd
{\bf Definition 2.} A topological conjugacy $\psi : X_\A \mapsto
X_{\A'}$ between tiling systems is a sliding block code if
for every $n'>0$ there is $n>0$ such that for every $x,\, y\in X_{\A}$
such that $B_{n}\cap \partial x=B_{n}\cap \partial y$ we have
$B_{n'}\cap \partial(\psi x)=B_{n'}\cap \partial(\psi y)$.
\vs.1 \nd
Our first result is:
\vs.1 \nd
{\bf Theorem 1.} Within the subclass of substitution tiling systems
with invertible substitution, every topological conjugacy is a sliding
block code.
\vs.1
Before defining the subclass of ``substitution'' tiling systems
in general we present some relevant examples.

A ``Penrose'' tiling of the plane, Figure 1, can be made as
follows. Consider the 4 (colored) tiles of Figure 2. Divide each tile
(also called a ``tile of level 0'') into 2 or 3 pieces as in Figure 2
and rescale by a linear factor of the golden mean
$\tau=(1+\sqrt{5})/2$ so that each piece is the same size as the
original.  This yields 4 collections of tiles that we call ``tiles of
level 1''.  Subdividing each of these tiles and rescaling gives 4
collections of tiles that we call tiles of level 2.  Repeating the
process $n$ times gives tiles of level $n$.  A Penrose tiling is a
tiling of the plane with the property that every finite subcollection
of tiles is congruent to a subset of a tile of some level.  A Penrose
tiling has only 4 types of tiles, each appearing in 10 different
orientations.

A ``pinwheel'' tiling of the plane, Figure 3, uses two basic tiles: a
1-2-$\sqrt{5}$ right triangle and its mirror image, as shown in
Figure 4 with their substitution rule [R1].
Notice that at the center of each tile of level 1 there is a tile of
level 0 similar to the level 1 tile but rotated by an angle $\alpha =
\tan^{-1}(1/2)$. Thus the center tile of a tile of level $n$ is
rotated by $n\alpha$ relative to the tile.  Since $\alpha$ is an
irrational multiple of $\pi$, we see, using the fact that within a
tile of level 2 there is a tile of level 0 similar and parallel to the level 2
tile, that this rotation never ends, and each tiling contains tiles
in infinitely many distinct orientations.

More generally, for any integers $m<n$ we consider the ``$(m,n)$-pinwheel''
tilings defined (for $m=3,\ n=4$) by the substitution of Figure 5, whose
tiles are $m$-$n$-$\sqrt{m^2+n^2}$ right triangles.  Like the ordinary
pinwheel, such variant pinwheel tilings also necessarily have tiles in
infinitely many distinct orientations.

It is easy to construct explicit examples of Penrose and pinwheel
tilings.  Pick a tile to include the origin of the plane.  Embed this
tile in a tile of level 1 (there are several ways to do this).  Embed
that tile of level 1 in a tile of level 2, embed that in a tile of
level 3, and so on.  The union of these tiles of all levels will cover
an infinite region, typically---though not necessarily---the entire
plane.

In order to generalize from the above examples we need some further
notation. A ``patch'' is a (finite or infinite) subset of a tiling
$x\in X(\A)$; the set of all finite patches for $\A$ will be denoted
$W_\A$.  A ``substitution function'' $\phi$ is a map from
$W_\A$ to itself defined by a decomposition of each
tile type, stretched linearly about the origin by a factor
$\lambda_\phi > 1$, into congruent copies of the original tiles. (Recall
the Penrose and pinwheel examples .) We assume:
\vs.1

\item{i)} For each $k>0$ there are only a finite number of possible
patches, up to Euclidean motion, obtained by taking a ball of radius $k$
around any point inside a tile $T$ of level $n$,
$\phi^n(T)$, where $T$ and $n$ are arbitrary. (This is usually called
``finite local complexity''.)

\item{ii)}  For each tile $T\in\A$, $\phi(T)$ contains at
least one tile of each type. (This is usually called ``primitivity''.)

\item{iii)} For every tile $T\in\A$ there is $n_T \ge 1$
such that $\phi^{n_T}(T)$ contains a tile of
the same type as $T$ and parallel to it.
\vs.1

Condition i) is highly significant, the remaining conditions much less so.
In interesting cases condition ii) can usually be obtained by
replacing the substitution by a power of itself; note that this does
not affect the tiling dynamical system at all.  Condition iii) is related
to the existence of a fixed point of the substitution; we know of no
interesting examples of systems not satisfying this condition.
\vs.1
\nd {\bf Definition 3.} For a given alphabet $\A$ of (possibly colored)
tiles, and a
substitution function $\phi$, the ``substitution tiling system'' is the
compact subspace $X_\phi\subset X(\A)$, invariant under $\G_E$, of
those tilings $x$ such that every finite subpatch of $x$ is congruent
to a subpatch of $\phi^n(T)$ for some $n>0$ and $T\in \A$. The
map $\phi$ extends naturally to a continuous map (again
denoted $\phi$) from $X_\phi$ into itself.
\vs.1

There are two natural relaxations of this definition of substitution
tiling systems. For any tiling system $X$ and any positive constant
$\lambda$, let $\lambda X$ denote the system of tilings obtained by
rescaling each tiling in $X$ by $\lambda$.  If $X_\phi$ is a
substitution tiling system, then $X_\phi$ and $\lambda_\phi X_\phi$
are topologically conjugate, via a sliding block code that
associates tiles (of level 0) of tilings in $\lambda_\phi X_\phi$ with
tiles of level 1 of tilings in $X_\phi$.

\vs.1
\nd {\bf Definition 4.} If a tiling system $X$ has the property that,
for some $\lambda>1$, $\lambda X$ and $X$ are topologically conjugate
via a sliding block code, then $X$ is a
``pseudo-substitution tiling system'', and the map $X \mapsto X$ obtained
by first rescaling by $\lambda$ and then applying the conjugacy is called
a ``pseudo-substitution''.  If $X$ and $\lambda X$ are topologically
conjugate (not necessarily via a sliding block code), then $X$ is
a ``quasi-substitution tiling system''.
\vfill \eject \nd
\hbox{} \vs-.05
\epsfig .55\hsize; 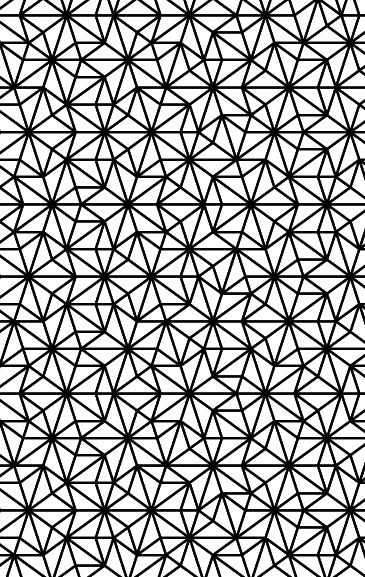;
\vs.2
\centerline{Figure 1. A Penrose tiling}
\vs.1
\hbox{}
\vbox{\epsfig .7\hsize; 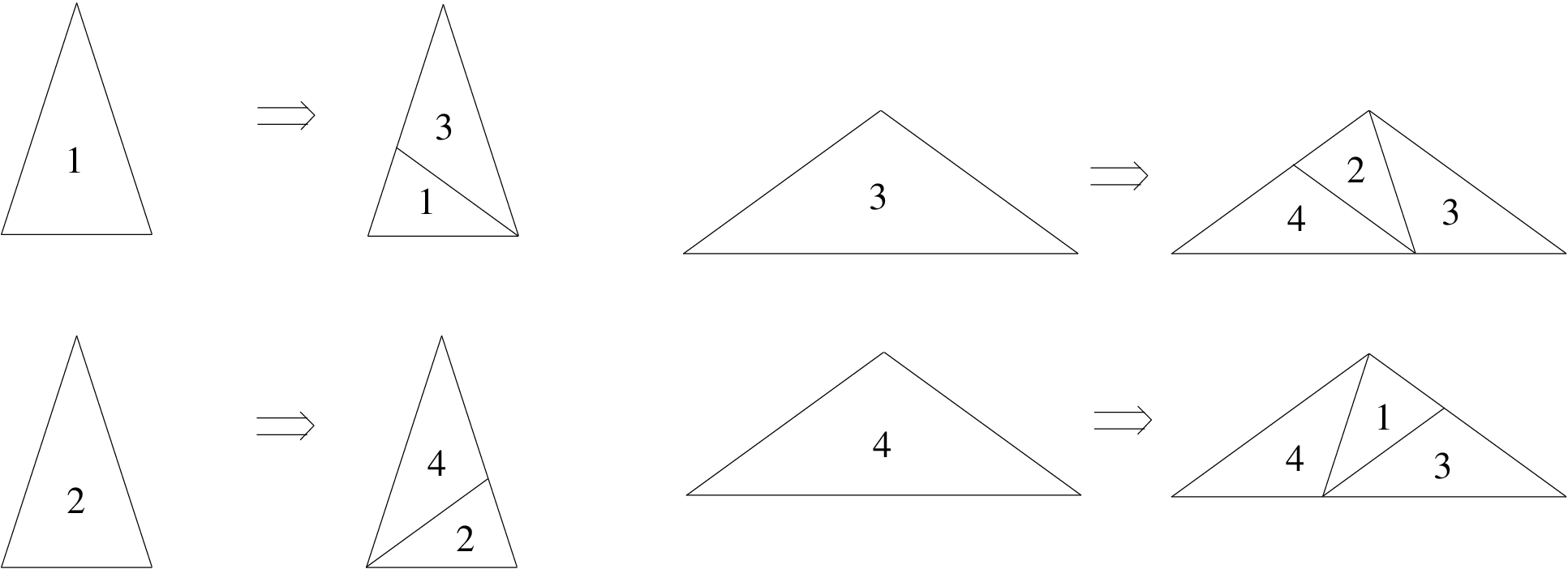;}
\vs.2
\centerline{Figure 2. The Penrose substitution}
\vfill\eject
\hbox{}\vs-.2
\hs.15 \epsfig .8\hsize; 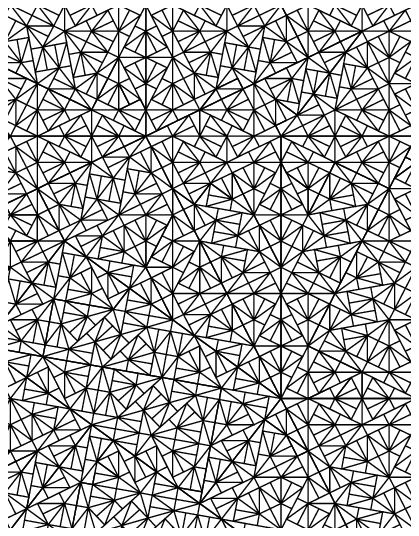;
\vs0 \centerline{Figure 3. A pinwheel tiling}
\hbox{}\vs.1
\vbox{\epsfig .5\hsize; 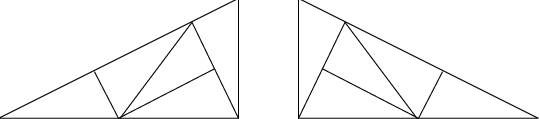 ;}
\vs.2
\centerline{Figure 4. The pinwheel substitution}
\vs1
\vbox{\epsfig .5\hsize; 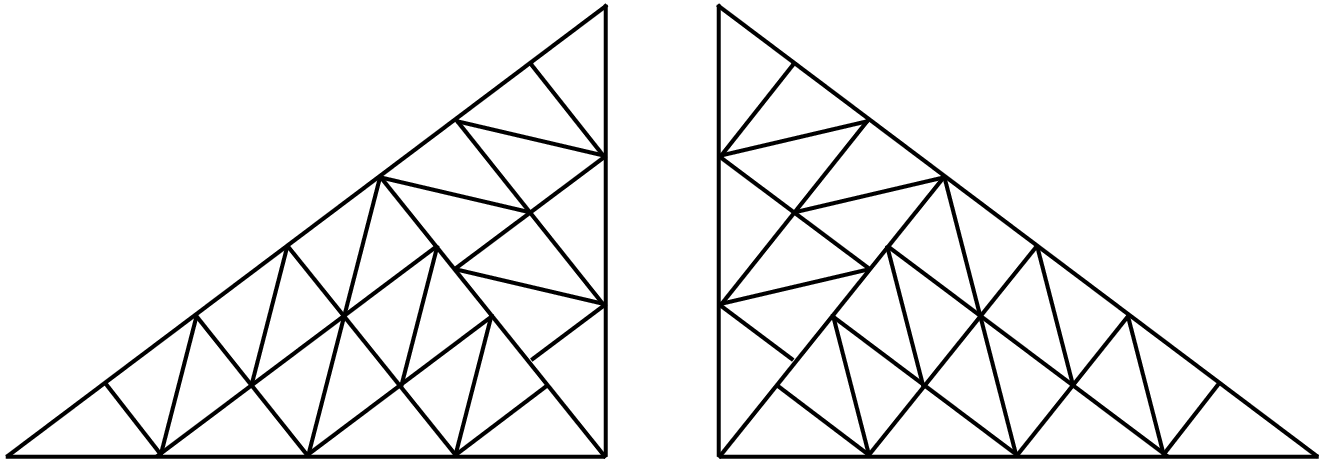 ;}
\vs.2
\centerline{Figure 5. The substitution for (3,4)-pinwheel tilings}
\vs.2
In the literature, pseudo-substitutions are sometimes called ``improper
substitutions'' or ``substitutions with amalgamation'', and their fixed
points are called ``pseudo-self-similar tilings''. In 2 dimensions, and
with some additional assumptions, the categories of substitution tiling
systems and of pseudo-substitution tiling systems are essentially
identical, thanks to a construction of Priebe and Solomyak [PS] that
converts a pseudo-self-similar tiling into a self-similar tiling. It is
generally believed that this construction can be generalized to higher
dimensions, dropping some of the restrictive assumptions, in which case
properties of substitution tiling systems, such as Theorem 1, can be
expected to apply to pseudo-substitution tiling systems.  However, the
following example (see also [P])
shows that the conclusions of Theorem 1 do {\it not}
apply to quasi-substitution tiling systems.

Following [RS2], we consider suspensions of the 1-dimensional
Fibonacci substitution subshift. That subshift is defined by the
alphabet $\B=\{0,1\}$ and the substitution of $0$ by $1$ and of 1 by
the word $0\,1$. One can make a family of suspensions of this subshift
by replacing 0 and 1 by marked closed intervals of any positive
lengths.  One gets a substitution tiling system if 0 is associated
with a segment $T_0$ of length $|T_0|= 1$ and 1 is associated with a
segment $T_1$ of length the golden mean $\tau=1+\sqrt{5}/2$,
$|T_1|=\tau$.  When $|T_1|/|T_0| \ne \tau$, the resulting tiling system
is merely a quasi-substitution tiling system [CS].
It was proven in [RS2] that two different Fibonacci
tiling systems, defined by $T_0, T_1$ and $T^\prime_0, T^\prime_1$,
are topologically conjugate if $|T_0|+\tau |T_1|=|T^\prime_0|+\tau
|T^\prime_1|$, but that {\it such a conjugacy cannot be a sliding
block code.}

Definitions 3 and 4 yield spaces of tilings in which the tiles
appear in all orientations, although tiles appear in at most countably
many different orientations within any {\it fixed} tiling.  Any
Penrose tiling, for example, has tiles in only 10 distinct orientations, but
the space of {\it all} Penrose tilings contains all rotated versions of any
tiling.

For each tiling $x$ and each $r>0$, consider the set of
Euclidean motions $g$ for which $x$ and $gx$ agree exactly on a ball of radius
$r$ around the origin. The subgroup of $SO(d)$ generated by the rotational
parts of the $g$'s is denoted $\G_\RO(r,x)$ and is called the
{\it relative orientation group} of $x$.  In [RS1] it was shown that
$\G_\RO(r,x)$ is independent of $r$ (and we henceforth write it as
$\G_\RO(x)$), and that the groups for different tilings
$x \in X_\phi$ are related by inner automorphism of $SO(d)$. There are no inner
automorphisms of $SO(2)$, so in 2 dimensions
the group is exactly the same for all tilings $x \in X_\phi$.

For the Penrose tiling, $\G_\RO= \Z_{10}$ is the set of rotations by
multiples of $2 \pi/10$. For the pinwheel tiling, $\G_\RO$ is generated
by rotations of $\pi/2$ and $2 \tan^{-1}(1/2)$.  For the $(m,n)$-pinwheel
tiling, the group is generated by rotations by $\pi/2$ and
$2 \tan^{-1} (m/n)$.

 From Theorem 1 we deduce that the relative orientation group is a conjugacy
invariant:
\vs .1
\nd {\bf Corollary 1.} If $\psi: X_1 \to X_2$ is a topological conjugacy
between substitution tiling systems with invertible substitution, and if
$x \in X_1$, then $\G_\RO(\psi(x))= \G_\RO(x)$.
\vs.1
As an application, consider the (1,2)-pinwheel and (3,4)-pinwheel
tilings. These have stretching factors $\sqrt 5$ and $5$, and
so cannot be distinguished by the homeomorphism invariant of
[ORS]. Moreover, their relative orientation groups are each
isomorphic to $\Z_4 \oplus \Z$ as abstract groups.  However, these
groups are different as subsets of $SO(2)$ (one is an index-2 subgroup
of the other), so the two tiling systems
cannot be topologically conjugate.

For each tiling $x \in X_\phi$, consider the closure of the orbit of
$x$ under translations.  Let $\G_{rel}(x)$ be the maximal subgroup of
$\G_E$ that maps this orbit closure to itself. It is not hard to see
that $\G_{rel}(x)$ is the smallest subgroup of $\G_E$ that contains
the closure of $\G_\RO(x)$ in $SO(d)$ and contains all translations.
If the translation orbit of $x$ is dense in $X_\phi$, as with the
pinwheel and its variants, then $\G_{rel}(x)=\G_E$.  We should add
that for tiling systems in which tiles only appear in finitely many
orientations in any tiling, it can be convenient to consider the
dynamical system of translations on the orbit closure of such a
tiling. Although our convention is to allow all of the Euclidean group
to act on tilings it is easy to obtain corresponding information about
these more limited dynamical systems from our results.

Our second  result concerns what is often called
the ``recognizability'' or ``unique
composition property'' of the substitution:
\vs.1 \nd
{\bf Theorem 2.} If for a substitution function $\phi$ there is some
tiling $x$ in $X_\phi$ not fixed
by any translation, for which the orbit under translation is dense in
$X_\phi$, then the extension of $\phi$ to $X_\phi$ is a homeomorphism.
\vs.1 \nd

Theorem 2 complements a result of Solomyak [S], which dealt with
tiling systems that had finite relative orientation groups.
Solomyak's result was itself a generalization of Moss\'e{}'s work on
1-dimensional subshifts [Mos].  For tilings in 3 or more dimensions
there is an additional case, where the relative orientation group is
infinite but not dense in $SO(d)$.  We have constructed a
pseudo-substitution tiling system in 3 dimensions, with $\G_{RO}$ a
dense subgroup of $SO(2)$, for which the pseudo-substitution is {\it not} a
homeomorphism. However, the recognizability of true substitutions in 3 or more
dimensions remains open.

Our third result concerns ``finite type''.  If $X$ is a tiling space and
$r>0$, let $X_r$ be the set of tilings for which every patch of radius $r$
also appears in some tiling in $X$.  If $r_1 > r_2$, then $X_{r_1}
\subseteq X_{r_2}$, and it is easy to show that $\cap_r X_r = X$.  If $X=X_r$
for some finite $r$, then we say that $X$ is of {\it finite
type}. Roughly speaking, this means that the patterns in $X$ are
defined by local conditions, whose range is at most $2r$.  For
subshifts, it is well known that being of finite type is an invariant
of topological conjugacy. (See [RS2] for an explicit proof of this
folk theorem).  We extend this to tiling systems:
\vs.1 \nd
{\bf Theorem 3.} Let $X,Y$ be topologically conjugate tiling systems,
each of finite local complexity.  $X$ is of finite type if and only if
$Y$ is of finite type.

\vs.2
\nd
{\bf II. Proofs of Theorems 1 and 2, and some related results}
\vs.1
We begin with the proof of Theorem 2.  We abbreviate $X_\phi$ by
$X$ and call a patch admissible if there is a tiling $ x\in X$
containing it. We assume $\phi$ has a fixed point in $ X$, not fixed
by any translation, whose orbit under translations is dense in $X.$
(The existence of a periodic point for $\phi$ follows from iii), and
we are free to replace $\phi$ with a higher power.)  Let $
H:\G_E\mapsto\G_E$ be defined by $\phi(gx)= H(g)\phi(x).$

The following four lemmas are proved with standard arguments, as sketched
below.
\vs.1 \nd
{\bf Lemma 1.} The extension $\phi: X\mapsto X$ is surjective.
\vs.1 \nd
Sketch of proof. Since
$X$ is compact, $\phi(X)$ is a closed subset of $X$. To see
that it is dense, note that any admissible patch is a subset of some tile of
level $n$. Thus, for any tiling $x \in X$, and any $r>0$, $B_r^x=B_r^{\phi(y)}$
for some tiling $y$. \qed
\vs.1 \nd
{\bf Lemma 2.} There is a constant $ C>0$ such
that for every $r>0$ and for every pair of admissible patches
$P,P'$ with $\hbox{ supp}(P)\subset B_r$ and
$B_{ C r}\subset\hbox{ supp}(P')$ there exists $g\in\G_E$ such that
$gP\subset P'.$
\vs.1 \nd Sketch of proof.
By finite local complexity and primitivity, there exists $N$ such that for
every $x\in X$ and for every prototile $T,$ the patch
$\phi^N(T)$ contains a congruent copy of $B_m^x.$
Take $C$ greater than $\lambda^{N+1}/m$ times the maximum diameter of a tile.
Suppose $r,P,P'$ are as in the statement of the lemma.  Let $n$ be the least
integer such that $r\lambda^{-n}\leq m.$  Then $P'$ contains a tile of 
level $N+n$, and every tile of level $N+n$ contains a congruent copy
of $P.$ \qed
\vs.1 \nd
{\bf Lemma 3.}
If $\b_1,\dots,\b_k\in\E^d$ and $\t=\sum_{i=1}^ka_i\b_i$ with
$a_i\in\N$ then there exist $\t_0,\t_1,\dots,\t_\ell$ such that
$\t_0=0,$ $\t_\ell=\t,$ and for each $j=1,\dots,\ell,$
$$
\t_j-\t_{j-1}\in\{\pm \b_i:i=1,\dots,k\} \cheqno\es
$$ \nd
and $\t_j$ lies within $(k/2)\max_{1\leq i\leq k}\|\b_i\|$ of the
straight-line path from {\bf O} to $\t.$
\vs.1 \nd
Sketch of proof. Each point along the straight-line path from {\bf O} to $\t$ is a linear
combination of the $\b_i$'s with real coefficients. Round these coefficients
to the nearest integer to get the sequence of $\t_j$'s.\qed
\vs.1 \nd
{\bf Lemma 4.} If $\G_{rel}(y)=\G_E$ there is a
constant $ D$ such that if $P,P'$ are
admissible patches in $y$ with $\hbox{ supp}(P)\subset B_r$ and
$B_{ D r}\subset\hbox{ supp}(P')$ then there exist
$\alpha_1,\dots,\alpha_n\in SO(d)$ and $\t_1,\dots,\t_n\in\R^d$ such that
$\alpha_iP+\t_i\subset P',$ $i=1,2,\dots,n,$ and such that no proper
subspace of $\E^d$ is invariant under all the $\alpha_i.$
\vs.1
\nd Sketch of proof. This is similar to Lemma 2.  One can prove it first
when $P$ is a tile and then extend the result to larger patches by
inverting the substitution.\qed
\vs.15
Let $ M=\max\{\hbox{diam}( T): T\in\A\}.$  Let us say that a patch $P$
has {\it period} $g\in\G_E$ if $P\cup gP$ is a patch, i.e.,
if $P$ and $gP$ agree  where their supports overlap (we do not
require that they actually overlap).  Alternatively, $P$ has period
$g\in\G_E$ if and only if whenever $ T\in P$ is such that
$g T^\circ$ intersects $\hbox{ supp}(P)$ we also have $g T\in P.$
Of course any subpatch of a patch of period $g$ has period $g.$
\vs.1 \nd
{\bf Lemma 5.}
If $\{\b_1,\dots,\b_{k-1}\}\subset\E^d$ is a basis for a lattice $\L$
and $P$ is a patch having all periods in $\L$ and additional
translational period $\b_k$ such that $B_r\subset\hbox{ supp}(P)$ where
$$
r>(d+1)\max_{1\leq i\leq k}\|\b_i\|+4 M \cheqno\et
$$ \nd
then $\langle \b_1,\dots,\b_k\rangle$ is a lattice of periods for
$B_{r/2}^P.$
\vs.1 \nd
Proof. Suppose $ T\in B_{r/2}^P$ and
$ T^\circ+\t$ intersects $\hbox{ supp}(B_{r/2}^P),$
where $\t=a_1\b_1+\dots+a_k\b_k$ with each $a_i\in\N.$  Let
$\b\in T\cap B_{r/2}.$  Then $\b+\t\in B_{r/2+2 M}.$  Let
$\t_0,\dots,\t_\ell$ be as in Lemma 3.  The
straight-line path from $\b$ to $\b+\t$ lies in $B_{r/2+2 M}$ so
each $\b+\t_j$ is in  $B_r.$  Thus, for each $j$ we have
$( T^\circ+\t_j)\cap P\neq\emptyset$ and by finite induction
$ T+\t_j\in P$.  It follows that
$ T+\t\in B_{r/2}^P.$\qed
\vs.1 \nd
{\bf Lemma 6.}
If $k<d$ and $\{\b_1,\dots,\b_k\}$ is a basis for a lattice $\L$
and $P$ is an admissible patch having all periods in $\L$ such
that $B_r\subset\hbox{ supp}(P)$ where
\vs-.05
$$
r> D((d+1)\max_{1\leq i\leq k}\|\b_i\|+4 M) \cheqno\eu
$$ 
\vs-.05 \nd
then there exists $\b_{k+1}\in\E^d\setminus\hbox{ span}(\L)$ with
$\|\b_{k+1}\|\leq\max_{1\leq i\leq k}\|\b_i\|$ and such that
$\langle \b_1,\dots,\b_{k+1}\rangle$ is a lattice of periods for
$B_{{r}/{3 D}}^P.$
\vs.1 \nd
Proof. Let $P'=B_{r/ D- M}^P.$  By Lemma 4,
there exist $\alpha\in SO(d)$ and $\t\in\R^d$ such that
$\alpha P'+\t\subset P$ and $\alpha\hbox{ span}(\L)\neq\hbox{ span}(\L).$
It follows that $P'$ has all periods in $\L$ as well as all periods
in $\alpha^{-1}\L.$  Let
$\b_{k+1}\in\alpha^{-1}\{\b_1,\dots,\b_k\}\setminus\hbox{ span}(\L).$  By
Lemma 5, $\langle \b_1,\dots,\b_{k+1}\rangle$ is a lattice of
periods for $B_{{r}/{2 D}- M/2}^P.$
Since ${r}/{2 D}- M/2>{r}/{3 D},$ the proof is complete.\qed
\vs.1
We next show that a large patch cannot have periods which are small
relative to the size of the patch.
\vs.1 \nd
{\bf Proposition 1.}
There is a constant $ K>0$ such that if $P$ is an
admissible patch whose support contains a ball of radius $r$
then every non-identity period $g$ of $P$ satisfies
$\|g\b-\b\|> K r$ for some $\b\in\hbox{ supp}(P).$
\vs.1 \nd
Proof. It suffices to prove there is a $K$ which satisfies the
conclusion for all sufficiently large $r.$ Recall the notation of $m$
as the inner radius. No tiling in $ X$ has a translational period of
magnitude less than $ m.$ Let $ K>0$ be less than each of the
following:

\item{a)}$4 C\lambda(\lambda+1)
(3 D)^{d-1}(d+4 M/ m)^{-1},$
\item{b)}$\inf\{\|\alpha-I\|_{\hbox{operator}}:
I\neq\alpha\in SO(d)\hbox{ and }\alpha\hbox{ fixes some element of } X\},$
\item{c)}${1\over 4 C}.$
\vs.1
Suppose there exist $r>4 M C$ and a patch $P$ having a
non-identity period $g\in\G_E$ with $\|g\b-\b\|\leq K r$ for all
$\b\in\hbox{ supp}(P),$ and such that $B_r\subset\hbox{ supp}(P).$  Let
$x\in X$ be a fixed point for $\phi,$ let
$P'=B_{r/ C- M}^x,$ and let $h\in\G_E$ be
such that $hP'\subset P.$  Then $h^{-1}gh$ is a period for $P'$
and
$$
\|h^{-1}gh\b-\b\|\leq K r, \ \hbox{ for all }\b\in P'. \cheqno\ev
$$
Since $\phi(P')\supset P'$ and $ H(h^{-1}gh)$ is a period
for $\phi(P')$, $ H(h^{-1}gh)(h^{-1}gh)^{-1}$ is a period
for $P'\cap(h^{-1}gh)^{-1}P'.$  We have
$$
B_{r({1\over C}- K)- M}^x\subset
P'\cap(h^{-1}gh)^{-1}P' \cheqno\ew
$$ \nd
and
$$
\| H(h^{-1}gh)(h^{-1}gh)^{-1}\b-\b\|\leq (\lambda+1) Kr,
\hbox{ for all }\b\in\hbox{ supp}(B_{r({1\over C}- K)- M}^x). \cheqno\ex
$$
Now $ H(h^{-1}gh)(h^{-1}gh)^{-1}$ is a translation, say by $\b_1\in\R^d,$
and if $\b_1=0$ then $h^{-1}gh\in SO(d)$ and $h^{-1}ghb=b.$ This last
is impossible due to our choice of $ K,$ hence
$$
0<\|x_1\|\leq(\lambda+1) K r. \cheqno\ey
$$

An application of Lemma 5 followed by $d-1$ applications
of Lemma 6 (we will see in a moment that $r$ is large enough
for this) yields a $d$-dimensional lattice
$\L=\langle \b_1,\dots,\b_d\rangle$ of periods for
$P'=B_{r'}^x$, where
$$
\|\b_i\|\leq(\lambda+1) K r,\qquad i=1,\dots,d, \cheqno\ez
$$ \nd
and
$$\eqalign{
r'&={1\over 2(3 D)^{d-1}}
r({1\over C}- K)- M\cr
&\geq{r\over 4 C(3 D)^{d-1}}\cr
&\geq(d+4 M/ m)\lambda(\lambda+1) K r\cr
&\geq(d+4 M/ m)\lambda\max_{1\leq i\leq d}\|x_i\|.\cr} \cheqno\eaa
$$
Thus $r'>\lambda(\sum_{i=1}^d\|x_i\|+ M),$ so the
fundamental domain $F=\{t_1\b_1+\dots+t_d \b_d:t_i\in [0,1]\,\}$
for $\L$ is such that
$\phi(F^x)\subset B_{r'}^x.$  This
implies that every tile in $\phi(F^x)$ is a translate by
an element of $\L$ of a tile in $F^x,$ and it follows
that all tiles in $ x$ are translates of tiles in
$F^x,$ contradicting that $x$ should have tiles in
infinitely many different orientations.\qed
\vs.1
For a tile $ T$ and $n\geq 0$ let $\P_n( T)$ be the set of
admissible patches $P$ for which $\phi^n( T)\subset\phi^n(P)$
and $\phi^n( T)\not\subset\phi^n(P')$ for any proper subpatch $P'$ of
$P.$  Then each $\P_n( T)$ is finite and
$\{ T\}=\P_0( T)\subset\P_1( T)\subset\cdots.$
\vs.1 \nd
{\bf Lemma 7.} For each tile $ T$ there is a positive integer $N_T$ such
that $\P_{N_T}( T)=\P_{N_T+1}( T)=\cdots.$

\vs.1 \nd
Proof. Set $\P( T)=\cup_{n\geq 0}\P_n( T)$ and let
$r>0,\ {\bf y}\in\E^d$ such that $B_r({\bf y})\subset T.$  By finite local
complexity, $\P( T)$ has only finitely many patches up to rigid motion,
since every $P\in\P( T)$ is of the form $({T^\circ})^x$
for some tiling $x\in X.$

If $P,gP\in\P_n( T)$ for some $g\in\G_E$ with
$\|g\b-\b\|< K r$ for all $\b\in P$ then $ H^n(g)$ is a period
for $\phi^n(P)$ which violates Proposition 1.  Thus
for each patch $P\in\P( T)$ the set of $g\in\G_E$ for which
$gP\in\P( T)$ is discrete and bounded, hence finite.  It follows that
$\P( T)$ is finite, which is equivalent to the desired result.\qed

\vs.1 \nd
{\bf Lemma 8.}
Suppose $\P_n( T)=\P_{n+1}( T).$  If $x\in X$ is such
that $\phi^{n+1}( T)\subset\phi(x)$ then
$\phi^n( T)\subset x.$
\vs.1 \nd
Proof. Let $x'\in X$ such that $\phi^n(x')=x.$
Then $\phi^{n+1}( T)\subset\phi(x)=\phi^{n+1}(x'),$
hence there exists $P\in\P_{n+1}( T)$ such that
$P\subset x'.$  Since $P\in\P_n( T)$ we have
$\phi^n( T)\subset\phi^n(P)\subset\phi^n(x')=x.$\qed

\vs.1
\noindent{\bf Proof of Theorem 2.}
Set $N=\max\{N_ T: T\in\A\}.$  Let $x\in X,$
and let $x_1, x_2\in X$ be any tilings such that
$$
\phi(x_1)=\phi^{N+1}(x_2)=x. \cheqno\eab
$$
We only need to show that $x_1=\phi^N(x_2)$.

Let $ T\in x_2$ be any tile and let $g\in\G_E$ such that $g T$
is a tile in $\A$.  Then
$$
\phi^{N+1}(g T)\subset\phi^{N+1}(gx_2)=\phi( H^N(g)x_1). \cheqno\eac
$$
By Lemma 8, $\phi^N(g T)\subset H^N(g)x_1,$
and hence $\phi^N( T)\subset x_1$.\qed

\vs.1 \nd
{\bf Remarks.}
The preceding arguments show that Proposition~1 is tantamount to
recognizability, even for substitutions that do not satisfy hypothesis iii)
and hence do not have a fixed point.  This equivalence will be used in 
the proof of Propositions 3 and 4, below.
The existence of a fixed point
and the fact that $\G_{rel}=\G_E$ were used to prove Proposition~1.
We believe the conclusion is false without the latter assumption.  
\vs.1

We now begin the proof of Theorem 1.
For $i=1,2,$ let $\phi_i$ be a substitution on alphabet $\A_i$ with linear
scaling factor $\lambda_i$ such that $\phi_i: X_{\phi_i} \mapsto
X_{\phi_i}$ is a homeomorphism.  Write $X_i$ for $X_{\phi_i}$ and suppose
$\psi:(X_1,\G_E)\mapsto (X_2,\G_E)$ is a topological conjugacy. 
\vs.1 \nd
\noindent{\bf Notation.} For $r>0$, $\a \in \E^d$ and a tiling $y$,
$B_r(\a)^y$ denotes the patch of $y$ consisting of all tiles in $y$
that intersect the open ball of radius $r$ about $\a$. We abbreviate
$B_r({\bf O})^y$ as $B_r^y$.  Pick an ``inner radius'' $m$ such that every
tile in $\A_1\cup\A_2$
contains an open ball of radius $m$.
For patch-valued functions $P,Q$ on $X_1$
we say $P$ determines $Q$ (or $Q$ is determined by $P$) if whenever
$x,y\in X_1$ and $P(x)=P(y)$ we also have $Q(x)=Q(y).$
\vs.1
It follows from the fact that $\phi_i:X_i\mapsto X_i$ is a homeomorphism
that there is a ``recognizability radius'' $D_i>0$ such that for $x\in X_i$ the
patch $B_{D_i}^x$ determines the patch consisting of tiles containing the origin
in $\phi^{-1}_i(x).$
\vs.1
\nd {\bf Lemma 9.}
There is a constant $\rho>0$ such that if $n\in\N$ and
$r>\lambda_2^n\rho$ then for $y\in X_2$ the patch
$\displaystyle B_{r/(2\lambda_2^n)}^{\phi_2^{-n}(y)}$
is determined by the patch $B_r^y.$
\vs.1 \nd
Proof. Take $\rho = 2D_2/(\lambda_2-1)$.
A patch of radius $r$ in  $y$ determines a
patch of radius $(r-D_2)/\lambda_2$ in $\phi_2^{-1}(y)$, of
radius $[(r-D_2)/\lambda_2 - D_2]/\lambda_2$ in $\phi_2^{-2}(y)$,
and $r/\lambda_2^n - D_2(\lambda_2^{-1}+ \lambda_2^{-2} + \cdots +
\lambda_2^{-n})$ in $\phi_2^{-n}(y)$.
This last radius is greater than $r\lambda_2^{-n} - {D_2}/{(\lambda_2-1)}$,
which in turn is at least ${r}/{2 \lambda_2^n}$.\qed

\vs.1
Any element $g\in\G_E$ can be written uniquely as the composition of a
rotation and a translation, i.e., there exist unique $\alpha\in SO(d)$ and
$\s\in\R^d$ such that
$$
g\a=\alpha \a+\s \hbox{ for all }\a\in \E^d \cheqno\ec
$$
and we set
$$
\ell(g)=\|\alpha-I\|_{\hbox{operator}}+\|\s\|. \cheqno\ed
$$
\vs.1
\noindent{\bf Notation.} For patch-valued functions $P,Q$
on $X_1$ the phrase
``$P$ determines $Q$ up
to motion by some $g\in\G_E$ with $\ell(g)\leq\eta$''
means that if $x,y\in X_1$ are such that
$P(x)=P(y)$ then there exists $g\in\G_E$ with
$\ell(g)\leq\eta$ such that $Q(x)=gQ(y).$
\vs.1 \nd
{\bf Lemma 10.} There exist a constant $S_0>0$ and
a function $\eta:\R_+\mapsto \R_+$ such that $\lim_{r\to\infty}\eta(r)=0$
and if $r>S_0$ then for $y\in X_1$
the patch $B_r^y$ determines
the patch $B_{r-S_0}^{\psi (y)}$ up
to motion by some $g\in\G_E$ with $\ell(g)\leq\eta(r).$
\vs.1 \nd
Proof.
By uniform  continuity of $\psi$, there is a radius $S_0$ such that the patch
$B_{S_0}^y$ determines the tile at the origin of $\psi(y)$ and its immediate
neighbors, up to motion by less than $m/2$.
Since $\psi$ is a conjugacy, for any point
$\a \in \E^d$, the patch $B_{S_0}(\a)^y$ determines the tile at $\a$ in
$\psi(y)$ and its nearest neighbors,
up to a small motion. Applying this to all points $\a \in B_{r-S_0}$, we
have that the patch $B_{r}^y$ determines $B_{r-S_0}^{\psi(y)}$
up to an overall small motion $g$.  The bound on
$\ell(g)$ follows from uniform continuity of $\psi$.\qed
\vs.1
For $n\in\N$ put
$$
k(n)=\left\lfloor n{\log\lambda_1\over \log\lambda_2}\right\rfloor \cheqno\ee
$$
and set
$$
\psi_n=\phi_2^{-k(n)}\circ\psi\circ\phi_1^n. \cheqno\ef
$$

Fix $x\in X_1.$
Note that $\displaystyle {\lambda_1^n/\lambda_2^{k(n)}} \in [1, \lambda_2)$,
and if $\log\lambda_1/\log\lambda_2$ is rational then the sequence
$(\lambda_1^n\lambda_2^{-k(n)})_{n=1}^\infty$ is periodic and
$1$ is in its range, while if $\log\lambda_1/\log\lambda_2$ is irrational then
the range of this sequence is dense in $[1,\lambda_2].$  Thus there exists
a subsequence $n_i$ such that $\psi_{n_i}(x)$ converges, say to $x'\in X_2,$
and such that
$$
\lambda_1^{n_i}\lambda_2^{-k(n_i)}\to 1\qquad\hbox{as }i\to\infty. \cheqno\eg
$$
{\bf Proposition 2.} The sequence $\{\psi_{n_i}\}_{i=1}^\infty$ converges uniformly
to a conjugacy $\psi':(X_1,\G_E)\mapsto (X_2,\G_E)$
such that, for $r>\rho + S_0$ and for $y\in X_1$ the patch
$B_r^y$ determines the patch
$B_{r/2}^{\psi'(y)}.$
\vs.1 \nd
Proof. Step 1. If $g\in\G_E$ then
$\psi_{n_i}(gx)\to gx'$ as $i\to\infty.$  Indeed, if we
write $gx=\alpha x+\s$ with $\alpha\in SO(d)$ and $\s\in\R^d$ then
$$\eqalign{
\psi_{n_i}(gx)
&=(\phi_2^{-k(n_i)}\circ\psi )(\phi_1^{n_i}(gx))\cr
&=(\phi_2^{-k(n_i)}\circ\psi )
        (\alpha\phi_1^{n_i}(x)+\lambda_1^{n_i}\s)\cr
&=\phi_2^{-k(n_i)}(
\alpha\psi[\phi_1^{n_i}(x)]+\lambda_1^{n_i}\s)\cr
&=\alpha\psi_{n_i}(x)+\lambda_2^{-k(n_i)}\lambda_1^{n_i}\s\cr} \cheqno\eh
$$
which clearly converges to $gx'$ as $i\to\infty.$
\vs.1
\noindent Step 2.  Suppose $g,g'\in\G_E$ and $r>\rho+S_0$
are such that $B_r^{gx}=B_r^{g'x}$. For each $i,$
$$
B_{\lambda_1^{n_i}r}^{\phi_1^{n_i}(gx)}
=B_{\lambda_1^{n_i}r}^{\phi_1^{n_i}(g'x)}. \cheqno\ei
$$ \nd
By Lemma 10, there exists $h_i\in\G_E$ with
$\ell(h_i)\leq\eta(\lambda_1^{n_i}r)$ such that
$$
B_{\lambda_1^{n_i}r-S_0}^{h_i\psi [\phi_1^{n_i}(gx)]}
=B_{\lambda_1^{n_i}r-S_0}
^{\psi [\phi_1^{n_i}(g'x)]}. \cheqno\ej
$$ \nd
By Lemma 9
$$
B_{(\lambda_1^{n_i}r-S_0)/2\lambda_2^{k(n_i)}}
^{\phi_2^{-k(n_i)}(h_i\psi[\phi_1^{n_i}(gx)])}
=
B_{(\lambda_1^{n_i}r-S_0)/2\lambda_2^{k(n_i)}}
^{\phi_2^{-k(n_i)}(\psi[\phi_1^{n_i}(g'x)])}.\cheqno\ek
$$

Let $\alpha,\beta_i\in SO(d)$ and $\s,\t_i\in\R^d$ be such that
$g\a=\alpha \a+\s$ and $h_i\a=\beta_i\a+\t_i$ for all $\a\in\E^d$.
We have
$$
\phi_2^{-k(n_i)}(h_i\psi[\phi_1^{n_i}(gx)])
=\beta_i\alpha\psi_{n_i}(x)
+\lambda_2^{-k(n_i)}\lambda_1^{n_i}\beta_i\s
+\lambda_2^{-k(n_i)}\t_i, \cheqno\el
$$ \nd
and this converges to $gx'$ as $i\to\infty.$
We know from Step 1 that
$$
\phi_2^{-k(n_i)}(\psi[\phi_1^{n_i}(g'x)])
\to g'x'\hbox{ as }i\to\infty. \cheqno\em
$$ \nd
Since
$$
{\lambda_1^{n_i}r-S_0\over 2\lambda_2^{k(n_i)}}
\to{r\over 2}\qquad\hbox{ as }i\to\infty, \cheqno\en
$$ \nd
it follows that
$$
B_{r/2}^{gx'}
=B_{r/2}^{g'x'}. \cheqno\eo
$$
\vs.1
\noindent Step 3.  Define $\psi':X_1\mapsto X_2$ as follows.  Given
$y\in X_1$ let $\{g_j\}_j$ be a sequence in $\G_E$ such
that $g_jx\mapsto y$ as $j\to\infty.$  There exist group
elements $h_j$ tending to the identity in $\G_E$ and positive real
numbers $r_j$ tending to infinity such that
$B_{r_j}^{h_jg_jx}=B_{r_j}^y$
for each $j.$  By Step 2, $(h_jg_jx')_j$ converges, hence so does
$(g_jx')_j,$  and we define
$\psi'(y)=\lim_{j\to\infty}g_jx.$
The existence of this limit ensures that $\psi'$ is well defined and
continuous.  We see from Step 1 that it is a conjugacy, and the sliding
block code follows from Step 2.
\vs.1
\noindent Step 4.  We still need to show that $\psi_{n_i}$ converges to
$\psi'$ uniformly on $X_1.$
For $y\in X_1$ and $r>0$ one can
find $g\in\G_E$ such that $gx$
agrees with $y$ on a ball of radius $r$ about the origin. By
linear repetitivity (see Lemma 2),
we can always choose $g$ such that $\ell(g)< Cr$ for some fixed
constant $C$.  By the triangle inequality,
$$\eqalign{
m[\psi_{n_i}(y),\psi'(y)]
\leq
m[\psi_{n_i}(y),&\psi_{n_i}(gx)]
+m[\psi_{n_i}(gx),g\psi_{n_i}(x)]\cr
&+m[g\psi_{n_i}(x),\psi'(gx)]
+m[\psi'(gx),\psi'(y)].\cr} \cheqno\ep
$$ \nd
Given $\epsilon$, we will show that for $i$ large (with estimates
independent of $y$), and for the correct
choice of $g$, each term on the right hand side is bounded by $\epsilon/4$.
The left hand side (which is independent of the choices made) is then
bounded by $\epsilon$. Since the estimates on $i$ were independent of $y$,
the left hand side goes to zero as $i \to \infty$ at a rate that is
independent of $y$.

The argument of Step 2 shows
that the maps $\psi_{n_i}$ are uniformly continuous with estimates that
are independent of $i$.
As a result, the first term can be made small,
independent of $i$ (and $y$), by choosing $r$ large
enough.  The last term is also
small if $r$ is large, since $\psi'$ is uniformly continuous.  For fixed
$r$ (and hence fixed $g$), the
second term is bounded by $Cr |1 - \lambda_1^{n_i} \lambda_2^{-k(n_i)}|$,
which is small if $i$ is large enough.
Finally, the third term is small once $i$ is big enough that
$\psi_{n_i}(x)$ and $\psi'(x)$ agree up to a small motion on
a ball of radius $\gg Cr$ about the origin.\qed
\vs.1 \nd
{\bf Proposition 3.} For all sufficiently large $i$ we have
$\lambda_1^{n_i}\lambda_2^{-k(n_i)}=1.$
\vs.1 \nd
{\bf Proposition 4.} There exists $I\in\N$ such that
for all $i\geq I$, for all $y\in X_1,$
\vs-.05
$$
\psi_{n_i}(y)=\psi'(y)+\s_{y,i} \cheqno\eq
$$
\vs-.05 \nd
for some $\s_{y,i}\in\R^d.$
\vs.1 \nd
Proof of Propositions 3 and 4.
For fixed $r,\varepsilon>0$ to be specified in the proof one can find
$\delta>0$ such that if
$y,y'\in X_1$ with $d(\psi_{n_i}(y),\psi'(y))<\delta$ then
there exists $g_{y,i}\in\G_E$ with $\ell(g_{y,i})<\varepsilon$
such that $\psi_{n_i}(y)$ and $g_{y,i}\psi'(y)$ agree on the ball
of radius $r$ centered at the origin.
Choose $I$ such that
$d(\psi_{n_i}(y),\psi'(y))<\delta$
for all $i\geq I,$ $y\in X_1.$ Let us consider $i\geq I.$
If $r$ is chosen large enough and $\varepsilon$ small enough,
then $g_{y,i}$ is uniquely determined by the above conditions and
varies continuously with $y,$ for otherwise
we would have a large patch with a small period, contradicting
the recognizability hypotheses (see the
remarks following the proof of 
Theorem 2).  Let $\alpha_{y,i}\in SO(d)$ and $\s_{y,i}\in\R^d$
denote the rotational and translational parts of $g_{y,i},$ respectively.
We have, for $\t\in\R^d$ with $\|\t\|<r/\lambda_2,$
\vs-.05
$$
\psi_{n_i}(y+\t)=\psi_{n_i}(y)+\lambda_1^{n_i}\lambda_2^{-k(n_i)}\t, \cheqno\eqq
$$
and this agrees with
$$
g_{y,i}\psi'(y)+\lambda_1^{n_i}\lambda_2^{-k(n_i)}\t
=g_{y,i}\psi'(y+\t)+\lambda_1^{n_i}\lambda_2^{-k(n_i)}\t-\alpha_{y,i}\t \cheqno\eqqq
$$
on the ball of radius $r-\lambda_2\t$ about the origin.
If $r$ is large enough and $\varepsilon$ small enough then this implies,
for all $\|\t\|$ sufficiently small, for all $y\in X_1,$
$$
\alpha_{y+\t,i}=\alpha_{y,i} \cheqno\eqqqq
$$
\vs-.05 \nd
and
$$
\s_{y+\t,i}=\s_{y,i}+\lambda_1^{n_i}\lambda_2^{-k(n_i)}\t-\alpha_{y,i}\t. \cheqno\eqqqqq
$$
By continuity $\alpha_{y+\t,i}=\alpha_{y,i}$ for all $\t,$ and this implies
that the above formula for $\s_{y+\t,i}$ holds for all $\t$ as well.  Now
$\s_{y+\t,i}<\varepsilon$ for all $\t,$ and this is only possible if
$\alpha_{y,i}$ is the identity and $\lambda_1^{n_i}\lambda_2^{-k(n_i)}=1.$
Thus $\psi_{n_i}(y+\t)=\psi'(y)+\s_{y,i}+\t.$
\qed
\vs.1 \nd
{\bf Corollary 2.} For each $i\geq I$ and each $y\in X_1$,
there exists $g_{y,i}$ in the center of
$\G_{rel}(y)$ such that $\psi_{n_i}(y)=g_{y,i}\psi'(y)$. Furthermore, if
$y'$ is in the closure of the translation orbit of $y$, then
$g_{y',i}=g_{y,i}.$
\vs.1
\nd Proof.
Fix $i$ and $y.$ By Proposition 4 there exists a translation $g_{y,i}$
such that $\psi_{n_i}(y)=g_{y,i}\psi'(y).$ Since $\psi_{n_i}$ and
$\psi'$ are conjugacies and all translations commute, we have
$g_{y,i}=g_{y',i}$ for any $y'$ in the translation orbit of $y.$ By
continuity, this last equality holds for all $y'$ in the closure of
the translation orbit of $y.$

To show that $g_{y,i}$ is in the center of $\G_{rel}(y),$ it suffices to show
that $\alpha g_{y,i}=g_{y,i}$ for every $\alpha\in SO(d)\cap\G_{rel}(y).$
Fix such $\alpha.$  By definition of $\G_{rel}(y)$ there is a sequence
$h_j\in\G_{rel}(y)$ such that $h_jy\to y$ as $j\to\infty$ and such that the
rotational part $\alpha_j$ of $h_j$ converges to $\alpha.$  If $g_{y,i}$ is
translation by $\s_{y,i}$ then
$$\eqalign{
\psi'(y)+\s_{y,i}
&=\psi_{n_i}(y)\cr
&=\lim_{j\to\infty}\psi_{n_i}(h_jy)\cr
&=\lim_{j\to\infty}h_j\psi_{n_i}(y)\cr
&=\lim_{j\to\infty}h_jg_{y,i}\psi'(y)\cr
&=\lim_{j\to\infty}h_j\psi'(y)+\alpha_j\s_{y,i}\cr
&=\lim_{j\to\infty}\psi'(h_jy)+\alpha_j\s_{y,i},\cr
&=\psi'(y)+\lim_{j\to\infty}\alpha_j\s_{y,i},\cr
} \cheqno\exxx
$$
and hence $\alpha_j\s_{y,i}\to\s_{y,i}.$  It follows that $g_{y,i}$ commutes
with $\alpha.$
\qed
\vs.1 \nd
{\bf Proposition 5.} If two tilings $y,y'$ in a substitution tiling space
$X_\phi$ agree on a single tile, then they are in the same translation
orbit closure.
\vs.1 \nd
Proof. Suppose that $y,y'\in X_\phi$ are tilings in different translation orbit
closures which agree on a tile; without loss of generality, we may assume that
the interior of the tile contains the origin, and thus,
$d(\phi^n(y),\phi^n(y'))\to0$ as $n\to\infty.$

There exists a tiling $z$ in the translation orbit
closure of $y$ and a rotation $\alpha\in SO(d)$ such that $\phi(z)=\alpha z.$
Let $\{n_i\}$ be an increasing sequence of integers such that 
$\alpha^{n_i}$ converges to the identity. 
Let $\beta\in SO(d)$ be such that $\beta z$ is in the translation orbit closure
of $y'.$  Then, for each $n_i\geq1,$ $\alpha^{n_i}y$ is in the 
translation orbit
closure of $\phi^{n_i}(y),$ and $\phi^{n_i}(\beta z)=\beta\alpha^{n_i}z$
is in the translation orbit closure of $\phi^{n_i}(y').$  Taking a limit
as $n_i \to \infty$, it follows that the distance from the translation
orbit closure of $y$ to that of $\beta z$ is zero, hence that $y$ and 
$\beta z$, and therefore $y'$, are in the same translation orbit closure.
\qed
\vs.1 \nd
{\bf Remark.} Using property iii), one can take $\alpha$ to be the
identity, thereby simplifying the proof of Proposition 5.  The above
argument, however, shows that the conclusions of Proposition 5 apply even
to substitutions that do not have a fixed point. 
\vs.1 \nd 
{\bf Proof of Theorem 1.}  
First we show that $\psi_{n_i}$ is a sliding block code
for $i \ge I$.  Suppose $x,y \in X_1$ agree on a large ball around the
origin, so that $\psi'(x)$ and $\psi'(y)$ agree on a (smaller) ball
around the origin. By Proposition 5, $x$ and $y$ lie in same translation
orbit closure.
However, $\psi_{n_1}$
and $\psi'$ differ by a (fixed) translation on this orbit closure, so
$\psi_{n_i}(x)$ and $\psi_{n_i}(y)$ agree on a (still smaller) ball around
the origin.

Now note that $\psi = \phi_2^{k(n_i)}
\psi_{n_i} \phi_1^{-n_i}$ is a
composition of three maps, each of which is a sliding block code
up to scaling. Each patch in $\psi(x)$ is determined by a
(much smaller) patch in $\psi_{n_i} \circ \phi_1^{-n_i}(x)$, which is
determined by a patch in $\phi_1^{-n_i}(x)$, which is determined by a
(much larger) patch in $x$.  Thus $\psi$ is a sliding block code.
\qed
\vs.1 \nd
{\bf Corollary 3.}
If the translation orbit of a tiling is dense in $X_1$,
then $\psi$ intertwines the actions of some powers of $\phi_1$ and $\phi_2.$
\vs.1
\nd Proof. In this case $\G_{rel}(y)$ of a tiling $y$ has no center, so $g_i$
is the identity and $\psi_{n_i}=\psi'$ for all $i \ge I$. But then
$$
\phi_2^{-k(n_i)}\circ\psi\circ\phi_1^{n_i}=\psi_{n_i}=\psi'=\psi_{n_{i+1}}
= \phi_2^{-k(n_{i+1})}\circ\psi\circ\phi_1^{n_{i+1}}. \cheqno\err
$$
Multiplying on the left by $\phi_2^{k(n_{i+1})}$ and on the right by
$\phi_1^{-n_i}$ gives
$$
\phi_2^{k(n_{i+1})- k(n_i)}\circ\psi =
\psi\circ\phi_1^{n_{i+1} - n_i}. \qed \cheqno\errr
$$
\vs.15 \nd
{\bf III. Proof of Theorem 3}
\vs.1
Suppose that $\psi:X \to Y$ is a topological conjugacy and that $Y$ is of
finite type.  By assumption, there exists a finite length ${r_1}$ such that
$Y_{r_1}=Y$.

As in the proof of Theorems 1 and 2, let $D$ be a length greater than the
diameter of any tile in either tiling system. By finite local complexity
there exists a radius $m>0$ such that
every tile contains a ball of radius $m$ and the only way to move two
adjacent tiles a distance $m$ or less and obtain an admissible local
pattern is to move the pair by a rigid motion.
Since $\psi^{-1}$ is
uniformly continuous, there is a length $r_2$ such that, for each $y \in Y$,
$B_{r_2}^y$ determines $B_{2D}^{\psi^{-1}(y)}$, up to a Euclidean
motion that moves each point in $B_{2D}$ by a distance of $m/2$ or less.
Likewise, there is a length $r_3$ such that, for each $x\in X$, $B_{r_3}^x$
determines $B_{r_1+r_2+3D}^{\psi(x)}$ up to a wiggle of size at most $m/2$.

We claim that $X_{r_3}=X$, and thus that $X$ is of finite type.
For if $x \in X_{r_3}$, then every patch of radius $r_3$ in $x$ corresponds
to an admissible patch in a tiling in $X$, and so determines a patch of
a tiling in $Y$ (up to a small motion).
That is, the tiling $x$ determines a combinatorial tiling
of the tiles of the $Y$ system, such that each patch of radius $r_1+r_2+3D$
is actually admissible. This local information can be stitched
together to form an actual tiling $y \in Y_{r_1+r_2+3D}=Y$.
The tiling
$\psi^{-1}(y)$ is then a tiling in $X$.  However, the combinatorial
structures of $x$ and $\psi^{-1}(y)$ are the same, since $B_{r_3}^x({\bf a})$
determines $B_{r_1+r_2+3D}^y$ up to a small rigid motion,
which determines $B_{2D}^{\psi^{-1}(y)}$ up to small rigid motion. Since
the tiles are rigid, this implies that $x$ is obtained by applying a
rigid motion to $\psi^{-1}(y)$, and is thus in $X$.  \qed
\vs.1 \nd
{\bf IV. Conclusions}
\vs.1
We have been concerned with topological conjugacy between tiling
dynamical systems, emphasizing the geometric aspects by including
systems in which the tiles appear, in each tiling, in infinitely many
orientations, thus incorporating the rotation group in an essential
way.  Some of our results are restricted to a subclass of tiling dynamical
systems, substitution systems, which can be thought of as
incorporating an extra group action which represents a certain
similarity: that is, not only have we extended the usual action of
the translation group by the rotations, we have in fact extended
further, to a subgroup of the conformal group.

Our first result is to show that topological conjugacies between
substitution systems with invertible substitutions are quite rigid.
We show that any conjugacy for the Euclidean actions automatically
extends to (some powers of) the similarities, and can be represented
by the natural analogue of a sliding block code. Tiling dynamical
systems are a geometric extension of subshifts, and these results all
have significant geometric meaning.

Our second result is that substitutions whose systems do not
admit periodic tilings are recognizable, as long as the relative
orientation group is either finite [S] or dense in
$SO(d)$. In particular, all nonperiodic substitutions in two dimensions are
recognizable. This result can then be used to generalize constructions such
as those in [PS] from the category of translationally finite tilings to the
more general case where tiles can appear in arbitrary orientation.

Part of the significance of
substitution subshifts and substitution tiling systems for dimension
$d\ge 2$ arises from the fact that, quite generally (see [Moz], [G]),
such a system carries a unique invariant measure and is measurably
conjugate to some uniquely ergodic system of finite type. Actually, the
proofs show more than measurable conjugacy; they show that off sets of
measure zero the map is bicontinuous, and it is boundedly finite to
one on the sets of measure zero. These associated finite type systems
are also of geometric interest, as part of the general effort of
understanding the symmetries of densest packings of bodies [R2].
Our third result is a step in this direction, showing that finite type
is a topological property among tilings with finite local complexity,
and not merely an artifact of the way one defines the tiles.

It would be significant if Theorem 1, and the conjugacy
invariants (the relative orientation groups) which follow from it,
apply to such tiling dynamical systems. For instance, we noted above
that the (1,2)-pinwheel and (3,4)-pinwheel systems cannot be
topologically conjugate.
We thus conclude with an open problem.
\vs.1 \nd
{\bf Open Problem.} Are the two finite type tiling systems, which are
measurably conjugate to the (1,2)-pinwheel and (3,4)-pinwheel systems,
topologically conjugate?

\vs.2 \nd
{\bf Acknowledgement} The authors are grateful for the support of
the Banff International Research Station, at which we formulated
and proved Theorem 3.

\vs.2
\nd {\bf V. References}
\vs.1


\item{[BSJ]} M. Baake, M. Schlottmann and P. D. Jarvis,
Quasiperiodic tilings with tenfold symmetry and equivalence with
respect to local derivability. {\it J. Physics A} {\bf 24} (1991),
4637-4654.

\item{[CS]} A. Clark and L. Sadun, When size matters: subshifts
and their related tiling spaces. {\it Ergodic Theory }${\&}$
{\it Dynamical Systems}, to appear.

\item{[G]} C.\ Goodman-Strauss, Matching rules and substitution tilings,
{\it Annals of Math.}, {\bf 147} (1998), 181-223.


\item{[LM]} D.\ Lind and B.\ Marcus, {\it An Introduction to Symbolic
Dynamics and Coding}, Cambridge University Press, Cambridge, 1995.

\item{[Mos]} B.\ Moss\'e, Puissances de mots et reconnaisabilit\'e des
point fixes d'une substitution, {\it Theor.\ Comp.\ Sci.} {\bf 99}
no.\ 2 (1992) 327-334.

\item{[Moz]} S.\ Mozes, Tilings, substitution systems and dynamical systems
generated by them, {\it J. d'Analyse Math.} {\bf 53} (1989), 139-186.

\item{[ORS]}N. Ormes,  C.\ Radin, and L.\ Sadun,
A homeomorphism invariant for substitution tiling spaces,
{\it Geometriae Dedicata} {\bf 90} (2002), 153-182.

\item{[P]} K.\ Petersen, Factor maps between tiling dynamical systems,
{\it Forum Math.} {\bf 11} (1999) 503-512.

\item{[PS]}  N.~Priebe and B.~Solomyak, Characterization of planar
pseudo-self-similar tilings. {\it Discrete Comput. Geom.} {\bf 26} (2001),
289--306.

\item{[R1]} C.\ Radin, The pinwheel tilings of the plane, {\it Annals of Math.} {\bf 139}
(1994), 661-702.

\item{[R2]} C.\ Radin,  Orbits of orbs: sphere packing meets Penrose tilings,
Amer. Math. Monthly (to appear).

\item{[RS1]} C.\ Radin and L.\ Sadun, An algebraic invariant for substitution tiling systems,
{\it Geometriae Dedicata} {\bf 73} (1998), 21-37.

\item{[RS2]} C.\ Radin and L.\ Sadun, Isomorphism of hierarchical structures,
{\it Ergodic Theory Dynam. Systems} {\bf 21} (2001), 1239-1248.

\item{[RW]} C.\ Radin and M.\ Wolff,
Space tilings and local isomorphism, {\it Geometriae Dedicata} {\bf 42} (1992), 355-360.

\item{[S]} B.\ Solomyak, Nonperiodicity implies unique composition for self-similar
translationally finite tilings, {\it Discrete Comput.\ Geom.} {\bf 20} (1998),
265-279.
\end